\documentstyle[12pt, amstex]{amsart}
\renewcommand{\pf   }{\noindent{\em Proof: }}
\newtheorem{Thm}{Theorem}[section]
\newtheorem{Def}[Thm]{Definition} \newtheorem{Rem}[Thm]{Remark}
\newtheorem{Lem}[Thm]{Lemma} \newtheorem{Cor}[Thm]{Corollary}
\newtheorem{Prop}[Thm]{Proposition}
\numberwithin{equation}{section}

\begin{document}

\fontsize{.5cm}{.5cm}\selectfont\sf

\vskip 0.5cm
\medskip

\title{On dual canonical bases}
\author{ Hechun Zhang }
\address{
Department of Mathematical Sciences, Tsinghua University,
  Beijing, 100084, P. R. China} \email{
hzhang@@math.tsinghua.edu.cn}

\footnotetext[1]{partially supported by NSF of China and by NWO of
the Netherlands, project number 613.006.573}

\begin{abstract}The
dual basis of the canonical basis of the modified quantized
enveloping algebra is studied, in particular for type $A$. The
construction of a basis for the coordinate algebra of the $n\times
n$ quantum matrices is appropriate for the study the
multiplicative property. It is shown that this basis is invariant
under multiplication by certain quantum minors including the
quantum determinant. Then a basis of quantum $SL(n)$ is obtained
by setting the quantum determinant to one. This basis turns out to
be equivalent to the dual canonical basis.
\end{abstract}

\maketitle

\medskip

\section{introduction}

\medskip

Throughout the paper, the base field is $K={\mathbb Q}(q)$, i.e.,
the field of quotients of polynomials in the indeterminate $q$
with rational coefficients. Let $A$ be an algebra over $K$. Two
elements $b, b^\prime\in A$ are called {\em equivalent} (denoted
by $b\sim b^\prime$) if there exists $m\in{\mathbb Z}$ such that
$b^\prime=q^m b$. Two elements $b, b^\prime$ are called {\em
$q$-commuting} if $bb^\prime\sim b^\prime b$.

Let $ g$ be the Kac-Moody algebra associated to a $n\times n$
symmetrizable Cartan matrix $A$. Let $U_q(g)$ be the quantized
enveloping algebra associated to $g$, with its two usual
subalgebras $U_q(n^+)$ and $U_q(n^-)$ (see section 2 for details).
The dual basis of the canonical basis of $U_q(n^-)$ has been
widely studied in literature. In \cite{k}, a conjecture posed by
Berenstein and Zelevinsky is stated as follows: Two elements $b_1,
b_2$ of the dual canonical basis are $q$-commuting with each
other, if and only if $b_1b_2\sim b$ for some $b$ in the dual
canonical basis. This property of basis is called the  {\em
multiplicative property}. By use of the Hall algebra technique,
the multiplicative property of the dual canonical basis of
$U_q(n^+)$ is studied in \cite{r}. In \cite{lec}, counter-examples
are given for the Berenstein-Zelevinsky conjecture by finding some
so-called imaginary vectors.  There are many connections between
the irreducible representations of  Hecke algebras of $A$ type and
the multiplicative property of the dual canonical basis, see
\cite{lec} and  \cite{lnt}.

Let $L(\lambda)$ be an irreducible highest weight module for
$U_q(g)$ and let $L^*(\lambda)$ be its graded dual. In \cite{lu1},
Lusztig constructed a canonical basis of the tensor product
$U(\lambda, \mu):=L(\lambda)\otimes L^*(\mu)$ which can be lifted
to a canonical basis $\tilde{B}$ of the so-called modified
quantized enveloping algebra $\tilde{U_q}(g)$. In the present
paper we will show that the module $L(\lambda)\otimes L^*(\mu)$ is
absolutely indecomposable if the Kac-Moody algebra $g$ is of
affine or indefinite type. Next, we focus on the case of type $A$.
By constructing a  basis of the coordinate algebra $O_q(M(n))$ of
the $n\times n$ quantum matrices, we get a basis of $O_q(SL(n))$
which turns out to be equivalent to the dual canonical basis. A
pleasant aspect of this construction is that it is appropriate to
study the multiplicative property of the basis.

{\bf Acknowledgement}\quad The major part of the present work was
done during the author visited University of Amsterdam. The author
would like to thank Professor T. Koornwinder for his hospitality.
The author  wishes to thank J. Du, H. P. Jakobsen, E. M. Opdam and
J. V. Stokman for valuable discussions.

\section{Kashiwara's construction}

\medskip

 Let $ g$ be the Kac-Moody algebra
associated to a $n\times n$ symmetrizable Cartan matrix $A$. One
can choose a bilinear form such that the integral weight lattice
is an even integral lattice. Let
$\Pi=\{\alpha_1,\alpha_2,\cdots,\alpha_n\}$ and
$\Pi^v=\{\alpha_1^v,\alpha_2^v,\cdots,\alpha_n^v\}$ be the set of
simple roots and the set of simple coroots respectively.  Let
$U_q( g)$ be the quantized enveloping algebra associated to $ g$
with generators $E_1, \cdots, E_n$, $F_1,\cdots, F_n$, $K_1,
K_1^{-1},\cdots, K_n, K_n^{-1}$ with the usual defining relations
(see e.g. \cite{k1}) by replacing $q$ by $q^2$  because we do not
want to use the square root of $q$ later.

 Let $U_q(n^+)$ (resp. $U_q(n^-)$) be the subalgebra
generated by $E_1, \cdots, E_n$ (resp. $F_1,\cdots, F_n$).
 For any dominant weight $\lambda$, denote by $L(\lambda)$ the irreducible highest weight
 module over $U_q(g)$ with the highest weight $\lambda$. Denote by
 $L^*(\lambda)$
 the graded dual of $L(\lambda)$ which is an irreducible lowest weight module
 with the lowest weight $-\lambda$.
Let $-$ be the automorphism of the algebra $U_q(g)$ given by
$$\bar{q}=q^{-1}, \bar{E_i}=E_i, \bar{F_i}=F_i,
\bar{K_i}=K_i^{-1}$$ for all $i$.
 Let $v_\lambda$ (resp. $v_\mu^*$) be a highest weight
 vector of $L(\lambda)$ (resp. the lowest weight vector of $L^*(\mu)$).
 Denote also by $-$ the linear
automorphism of the module $L(\lambda)$ and of the module
$L^*(\mu)$ given by
$$\overline{pv_\lambda}=\bar{p}v_\lambda,
\overline{pv_\mu^*}=\bar{p}v_\mu^*$$ for $p\in U_q(g)$.

\begin{Rem}Although we use $-$ to denote several different automorphism
of various spaces. One can identify the meaning of the $-$ from
the context.\end{Rem}

 In \cite{lu1}, Lusztig constructed a canonical basis of the tensor product
 $U(\lambda, \mu):=L(\lambda)\otimes L^*(\mu)$ which can be lifted to a canonical basis
 $\tilde{B}$
  of the modified
 quantized enveloping algebra $\tilde{U_q}(g)$.

 We will not go in detail about the canonical basis of
 the module $U(\lambda,\mu)$. However, we would like to show one
 remarkable fact about the module $U(\lambda, \mu)$. It is known that if $g$ is of
 finite type,  $U(\lambda, \mu)$ is of finite dimensional and is
 indecomposable if and only if one of $\lambda$ and $\mu$ is zero.
 However, if $g$ is of affine or indefinite type, the situation
 changes dramatically.

 \begin{Thm} If $g$ is of affine or indefinite type, then
 $$End_{U_q(g)}U(\lambda, \mu)\cong {\mathbb Q}(q).$$
 Hence, $U(\lambda, \mu)$ is absolutely indecomposable.\end{Thm}

\pf Clearly, if $\lambda$ or $\mu$ is trivial, then
$U(\lambda,\mu)$ is a lowest weight module or a highest weight
module and the theorem holds. Hence, we may assume that both
$\lambda$ and $\mu$ are nontrivial.

 It is known that
 $U(\lambda, \mu)$ is a cyclic module and is generated by $v_\lambda\otimes
 v_\mu^*$.
For any $\psi\in End_{U_q(g)}U(\lambda, \mu)$, then $\psi
(v_\lambda\otimes v_\mu^*) =u(v_\lambda\otimes v_\mu^*)\in
U(\lambda,\mu)_{\lambda-\mu}$ for some $u\in U_q(g)$ which is of
weight zero. If $u(v_\lambda\otimes v_\mu^*)$ is not a multiple of
$v_\lambda\otimes v_\mu^*$, then
$$u(v_\lambda\otimes v_\mu^*)=s(v_\lambda\otimes
v_\mu^*)+\sum_iu_iv_\lambda\otimes w_iv_\mu^*$$ where
$s\in{\mathbb Q}(q), u_i\in U_q(n^-), w_i\in U_q(n^+)$, for all
$i$, and the set $\{w_iv_\mu^*\}_i$ is linearly independent.
Choose $w_kv_\lambda^*$ such that its weight is maximal among all
of the weights of $w_iv_\lambda^*$ for all $i$. Assume that
$u_kv_\lambda\in L(\lambda)_\Lambda$, where $\Lambda$ must be
smaller than $\lambda$.

1. If the Cartan matrix $A$ is of indefinite type, then there
exists $\alpha^v_i$ such that $<\lambda-\Lambda, \alpha_i^v><0$,
i.e. $<\lambda,\alpha^v_i><<\Lambda, \alpha^v_i>$ and so
$F_i^{<\lambda,\alpha^v_i>+1}u_kv_\lambda\ne 0$. However,
$F_i^{<\lambda,\alpha^v_i>+1}u(v_\lambda\otimes
v_\mu^*)=\psi(F_i^{<\lambda,\alpha^v_i>+1}(v_\lambda\otimes
v_\mu^*))=0$. On the other hand,
$F_i^{<\lambda,\alpha^v_i>+1}u(v_\lambda\otimes v_\mu^*)=\sum_{m,
j}c_{ij}^{(m)}F_i^{<\lambda,\alpha^v_i>+1-m}u_jv_\lambda\otimes
F_i^mw_jv_\mu^*$, where $c_{ij}^{(m)}\in{\mathbb Q}(q)$. One can
see easily that $c_{ik}^{(0)}=1$. Hence,
$F_i^{<\lambda,\alpha^v_i>+1}u_kv_\lambda=0$. Contradiction!

2. Now, we may assume that the Cartan matrix $A$ is of affine
type. If there exists $\alpha_i^v$ such that
$<\lambda-\Lambda,\alpha^v_i><0$, then we can prove in the same
way as above. If $<\lambda-\Lambda,\alpha^v_i>\ge 0$ for all $i$,
then we must have $<\lambda-\Lambda,\alpha^v_i>=0$ for all $i$. As
there exists $E_i$ such that $E_iu_kv_\lambda\ne 0$, we have again
that $F_i^{<\lambda,\alpha^v_i>+1}u_kv_\lambda\ne 0$.

Therefore, $u(v_\lambda\otimes v_\mu^*)$ is a multiple of
$v_\lambda\otimes v_\mu^*$ and so $\psi$ is a scalar endomorphism.
Moreover, $U(\lambda,\mu)$  is  absolutely indecomposable. \qed

 Let $U^{{\mathbb Z}}(g)$ be the integral form of the quantized enveloping algebra
 which is a ${\mathbb Z}[q, q^{-1}]$-subalgebra of the quantized enveloping algebra
 $U_q(g)$ generated by the divided powers $E_i^{(s)}:=E_i^s/[s]!,
 F_i^{(s)}:=F_i^s/[s]!, K_i, K_i^{-1}$ for all $i$. The quantum integer is definied by
$[s]=\frac{q^{2s}-q^{-2s}}{q^2-q^{-2}}$,
  (one may refer to \cite{jan}
 for terminologies and notations).

Denote by $U_q(g)^*$ the linear dual of the algebra $U_q(g)$.
Since $U_q(g)$ is a $U_q(g)$ bi-module, $U_q(g)^*$ has an induced
$U_q(g)$ bi-module structure. Let
\begin{eqnarray}A_q(g):&=&\{f\in
U_q(g)^*|\text{ there exists }l\ge 0\text{ such that }\\\nonumber
E_{i_1}\cdots E_{i_l}f&=& fF_{i_1}\cdots F_{i_l}=0\text{ for any
}i_1,\cdots,i_l \}.\end{eqnarray}

The quantum Peter-Weyl theorem was proved in \cite{k}
\begin{Thm}As $U_q(g)$ bi-modules
$$A_q(g)\cong \oplus_{\lambda\in P} L(\lambda)\otimes L^*(\lambda)$$
where $u\otimes v\in L(\lambda)\otimes L^*(\lambda)$ viewed as a linear function on
$U_q(g)$ as follows:
$$(u\otimes v)(p)=<up, v>,\text{ for }p\in U_q(g),$$
where $L(\lambda)$ and $L^*(\lambda)$ are viewed as right $U_q(g)$
module and left $U_q(g)$ module respectively.
\end{Thm}

Let $A$ be the subring of ${\mathbb Q}(q)$ consists of the
rational functions of $q$ which are regular at $q=0$. Let $ -$ be
the ring endmorphism of ${\mathbb Q}(q)$ sending $q$ to $q^{-1}$.

Let $M$ be an integral $U_q(g)$ module. Then
$$M=\oplus_{\lambda} F_i^{(n)}(KerE_i\cap M_\lambda).$$
We define the lower Kashiwara operators $e_i, f_i$ of $M$ by
$$f_i(F_i^{(n)}u)=F_i^{(n+1)}u\text{ and
}e_i(F_i^{(n)}u)=F_i^{(n-1)}u$$ for $u\in KerE_i\cap M_{\lambda}$.

\begin{Def} A pair $(L, B)$ is called a lower crystal base of $M$
if it satisfies the following conditions:
\begin{enumerate}

\item  $L$ is a free sub-$A$-module of $M$ such that $M\cong
{\mathbb Q}(q)\otimes_A L$.

\item $B$ is a base of the ${\mathbb Q}$-vector space $L/qL$.

\item $e_iL\subset L$ and $f_iL\subset L$ for any $i$.

\item $e_iB\subset B\cup \{0\}$ and $f_iB\subset B\cup \{0\}$.

\item $L=\oplus_{\lambda\in P} L_{\lambda}$ and
$B=\cup_{\lambda\in P} B_\lambda$, where $L_{\lambda}=L\cap
M_{\lambda}$, $B_\lambda=B\cap L_\lambda/qL_\lambda$.

\item For any $b, b^\prime\in B$, $b^\prime=f_ib$ if and only if
$b=e_ib^\prime$. \end{enumerate}\end{Def}

The upper Kashiwara operators $e_i^\prime$ and $f_i^\prime$ are
define as follows: for $u\in KerE_i\cap M_\lambda$ and $0\le n\le
<\lambda, \alpha^v>$,
$$e_i^\prime(F_i^{(n)}u)=\frac{[<\lambda,
\alpha^v>-n+1]}{[n]}F_i^{(n-1)}u, $$
and
$$f_i^\prime(F_i^{(n)}u)=\frac{[n+1]}{[<\lambda,
\alpha^v>-n]}F_i^{(n+1)}u. $$

We say that $(L, B)$ is an upper crystal base if $(L, B)$
satisfies the the conditions in the definition of lower crystal
base with $e_i^\prime, f_i^\prime $ instead of $e_i, f_i$.

For $\lambda\in P$, we define $\psi_M\in AutM$ by
$$\psi_M(u)=q^{-2(\lambda, \lambda)}u$$
for $u\in M_\lambda$. It is known that
$\psi_M^{-1}e_i^\prime\psi_M$ (resp.
$\psi_M^{-1}f_i^\prime\psi_M$) coincides with $e_i$ (resp. $f_i$)
on $L/qL$.

In \cite{ka}, Kashiwara proved that

\begin{Lem} $(L, B)$ is a lower crystal base if and only if
$\psi_M(L, B)$ is an upper crystal base.\end{Lem}

Let $\cal L$$(\lambda)$ be the upper crystal lattice  which is the
smallest $A$ submodule of $L(\lambda)$ containing $v_\lambda$ and
is stable under the action of upper Kashiwara operators.
Similarly, let $\cal L$$^*(\lambda)$ be the upper crystal lattice
 which is the smallest
$A$ submodule of $L^*(\lambda)$ containing $v_\lambda^*$ and is
stable under the action of upper Kashiwara operators.
 Set
$${\cal L}(A_q(g)):=\oplus_{\lambda\in P_+}{{\cal L}}(\lambda)
\otimes {{\cal L}}^*(\lambda).$$
Define that
$$<\bar{u}, p>=
\overline{<u,\bar{p}>},$$ then one can check that
$\overline{u\otimes v}=\bar{u}\otimes \bar{v}$ for $u\in
L(\lambda)$ and $v\in L^*(\lambda)$. Hence
$$\overline{{\cal L}(A_q(g))}=\oplus_{\lambda\in P_+}\overline{\cal L}(\lambda)
\otimes \overline{\cal L}^*(\lambda).$$

Let
$$A^{\mathbb Z}_q(g)=\{f\in A_q(g)|<f,U^{\mathbb Z}(g)>\subset
{\mathbb Z}[q,q^{-1}]\}$$

Let $u\otimes v\in L(\lambda)\otimes L^*(\lambda)$, where $u$
(resp. $v$) is a weight vector of weight $\lambda_l$ (resp.
$\lambda_r$). Then $u\otimes v$ is called a weight vector with
left weight $\lambda_l$ and right weight $\lambda_r$. An element
in $A_q(g)$ is called a refined weight vector if it is a linear
combination of the elements $u\otimes v$ with the same left and
right weights.

Let us recall the definition of balanced triple. Let $V$ be a
vector space over ${\mathbb Q}(q)$, a $B$-lattice of $V$ is a
$B$-submodule $M$ of $V$ such that $V\cong {\mathbb Q}(q)\otimes_B
M$. Let $V_{\mathbb Z}$ be a ${\mathbb Z}[q,q^{-1}]$-lattice of
$V$, $L$ an $A$-lattice of $V$, and $\overline{L}$ an
$\overline{A}$-lattice of $V$. In \cite{ka}, it was proved that

\begin{Lem} Set $E=V_{\mathbb Z}\cap L\cap \overline{L}$. Then the
following conditions are equivalent.

\begin{enumerate}
\item $E\longrightarrow V_{\mathbb Z}\cap L/V_{\mathbb Z}\cap qL$
is an isomorphism.

\item $E\longrightarrow V_{\mathbb Z}\cap \overline{L}/ V_{\mathbb
Z}\cap q^{-1}\overline{L}$ is an isomorphism.

\item $V_{\mathbb Z}\cap qL\oplus V_{\mathbb Z}\cap
\overline{L}\longrightarrow V_{\mathbb Z}$ is an isomorphism.

\item $A\otimes E\longrightarrow L$, $\overline{A}\otimes
E\longrightarrow \overline{L}$, ${\mathbb
Z}[q,q^{-1}]\otimes_{\mathbb Z} E\longrightarrow V_{\mathbb Z}$,
${\mathbb Q}(q)\otimes_{\mathbb Z} E\longrightarrow V$ are
isomorphisms.
\end{enumerate}
\end{Lem}

We call $(L, \overline{L},V_{\mathbb Z})$ balanced if these
equivalent conditions are satisfied. Let us denote by $G$ the
inverse of the isomorphism $E\longrightarrow V_{\mathbb Z}\cap
\overline{L}/ V_{\mathbb Z}\cap q^{-1}\overline{L}$. If $B$ is a
base of $ V_{\mathbb Z}\cap \overline{L}/ V_{\mathbb Z}\cap
q^{-1}\overline{L}$, then $\{G(b)| b\in B\}$ is a base of $V$.

 In \cite{k}, it was proved that $(A_q^{\mathbb Z}(g), {\cal
L}(A_q(g)), \overline{{\cal L}(A_q(g))})$ is a balanced triple.
Hence there is a $\mathbb Z$ basis  $B^\prime$ of $$(A_q^{\mathbb
Z}(g)\cap {\cal L}(A_q(g))\cap \overline{{\cal L}(A_q(g))}).$$

 In \cite{k1}, it was shown that $B^\prime$ is the
dual basis of the canonical basis of the modified enveloping
algebra $\tilde{U}_q(g)$ if $g$ if of finite type.

In the following, we always assume that $g$ is of finite type. We
fix a reduced expression in the longest element of the Weyl group.
Let $F_{\beta_1},F_{\beta_2},\cdots,F_{\beta_N}$  be the ordered
root vectors given defined according to the chosen reduced
expression of the longest element in the Weyl group, where  $N$ is
the length of the longest element in the Weyl group. For any
$I=(i_1,i_2,\cdots,i_N)\in{\mathbb Z}_+^N$, denote by $F^I$ the
monomial $F_{\beta_1}^{(i_1)}F_{\beta_2}^{(i_2)}\cdots
F_{\beta_N}^{(i_N)}$ which form a  PBW type basis  of the
subalgebra $U_q(n^-)$. The monomial $E^I$ is defined similarly
which form a PBW type basis of the subalgebra   $U_q(n^+)$.

Let $B^-$ and $B^+$ be the canonical basis of $U_q(n^-)$ and
$U_q(n^+)$ respectively. For any dominant weight $\lambda$, denote
by
$$B^-_\lambda=\{b\in B^-|bv_\lambda\ne 0\}$$
and
$$B^+_\lambda=\{b^\prime\in B^+|b^\prime v_\lambda^*\ne 0\}.$$
Note that each dual canonical basis element
 is a refined weight vector.  Hence, we only need to consider the
 homogeneous part ${\cal L}(\lambda)_\mu\otimes {\cal L^*}(\lambda)_\gamma
 \cap\overline{{\cal L}(\lambda)_\mu\otimes {\cal
 L^*}(\lambda)_\gamma}
 \cap A_q^{\mathbb Z}(g)$.
 
   It is well-known that any canonical basis
element $b$ in $B^-$ is of the form
$$b=F^I+\sum_{I^\prime}a_{I,I^\prime}F^{I^\prime}$$
where the coefficients $a_{I, I^\prime}\in q{\mathbb Z}[q]$ and
the element $b$ is $-$ invariant. $F^I$ is called the leading term
of $b$. The canonical basis elements in $B^+$ have the similar
form.

Let
$$C_\lambda^-=\{F^I|F^I\text{ is the leading term of an
element }b\in B^-_\lambda\}$$
and let
$$C_\lambda^+=\{E^I|E^I\text{ is the leading term of an element
}b^\prime\in B^+_\lambda\}.$$ Then $C_\lambda^-v_\lambda$ (resp.
$C_\lambda^+v_\lambda^*$) is the PBW-basis of $L(\lambda)$ (resp.
$L(\lambda)^*$) with an order given by the chosen reduced
expression of the longest element in the Weyl group.

We order the PBW type basis $\bigcup_{\lambda} C_\lambda^-\otimes
C_\lambda^+$ by lexicographic ordering.

\begin{Thm}\label{dual} The basis $B^\prime$ is characterized by the following
two conditions.
\begin{enumerate}
\item $b^\prime=F^Iv_\lambda\otimes E^{I^\prime}v_\lambda^*+\sum_{I_k,
I_K^\prime} a_{I, I^\prime}^{I_k,I_k^\prime}
F^{I_k}v_\lambda\otimes E^{I_k^\prime}v_\lambda^*$ where $a_{I,
I^\prime}^{I_k,I_k^\prime}\in q{\mathbb Z}[q]$ and $a_{I,
I^\prime}^{I_k,I_k^\prime}\ne 0$ only if $(I_k, I_k^\prime)\le (I,
I^\prime)$, for any $b^\prime\in B^\prime$.
\item $\overline{b^\prime}=b^\prime$.\end{enumerate}\end{Thm}

\pf Clearly, each element $bv_\lambda\otimes b^\prime v_\lambda^*$
satisfies the two conditions. The uniqueness can be proved in the
same way as in \cite{d1} .\qed

The following result was proved in \cite{k}.

\begin{Prop} Let $x$ and $y$ be refined weight vectors of weights
$(\lambda_l,\lambda_r)$ and $(\mu_l,\mu_r)$ respectively. Then
$$\overline{xy}=q^{2(\lambda_r,\mu_r)-2(\lambda_l,\mu_l)}\bar{y}\bar{x}.$$
\end{Prop}

By using the above Proposition, one can easily verify that

\begin{Lem} The mapping
\begin{eqnarray}\phi: A_q(g) &\longrightarrow & A_q(g),\\\nonumber
q &\mapsto & q^{-1},\\\nonumber
 u\otimes v &\mapsto &
q^{((\lambda_l,\lambda_l)-(\lambda_r,\lambda_r))}\bar{u}\otimes\bar{v}.
\end{eqnarray}
 if $u\otimes v$  is of the left weight $\lambda_l$
and right weight $\lambda_r$,
 extends to an algebra
anti-automorphism of the algebra $A_q(g)$ over $\mathbb Q$.
\end{Lem}

Let $b^\prime\in B^\prime$ with weights $(\lambda_l,\lambda_r)$.
Then the element
$b=q^{\frac{1}{2}((\lambda_l,\lambda_l)-(\lambda_r,\lambda_r))}b^\prime$
is invariant under the anti automorphism $\phi$. Let
$$L^*=\{b|b^\prime\in B^\prime\}.$$
Then $L^*$ is also a $\mathbb Z$$[q,q^{-1}]$ basis of $A^{\mathbb
Z}(g)$.

It is clearly that the multiplicative properties of $B^\prime$ and
$L^*$ are the same.

\begin{Prop}Let $b_1,b_2\in L^*$. Assume that $b_1b_2\sim b$
 for some $b\in L^*$. Then $b_1b_2\sim b_2b_1$.\end{Prop}

 \pf Assume that $b_1b_2=q^ab$ for some $a\in\mathbb Z$. Applying
 the anti automorphism $\phi$, we deduce that
$b_2b_1=q^{-a}b$. Hence, $b_1b_2=q^{2a}b_2b_1$.\qed

\medskip

\section{The construction of the basis of $O_q(M(n))$}

\medskip

The coordinate algebra  $O_{q}(M(n))$ of the quantum matrix is an
associative algebra, generated by elements
$Z_{ij},i,j=1,2,\cdots,n$, subject to the following defining
relations:
\begin{eqnarray}\label{relations}Z_{ij}Z_{ik}&=&q^2Z_{ik}Z_{ij} \text{ if } j<k,\\
 Z_{ij}Z_{kj}&=&q^2Z_{kj}Z_{ij} \text{ if }i<k,\\
Z_{ij}Z_{st}&=&Z_{st}Z_{ij}\text{ if } i>s, j<t,\\
Z_{ij}Z_{st}&=&Z_{st}Z_{ij}+(q^2-q^{-2})Z_{it}Z_{sj} \text{ if } i<s, j<t.\end{eqnarray}
For any matrix $A=(a_{ij})_{i,j=1}^n\in M_n({\mathbb Z}_+)$ ( ${\mathbb Z}_+=\{0,1,\cdots\}$) we define a monomial $Z^A$  by
\begin{equation} Z^A=\Pi_{i,j=1}^nZ_{ij}^{a_{ij}},\end{equation}
where the factors are arranged in the lexicographic order on
$I(n)=\{(i,j)\mid i,j=1,\dots,n\}$.
It is well known that the set $\{Z^A|A\in M_n({\mathbb Z}_+\}$ is a basis of the algebra
$O_q(M(n))$.

From the defining relations (\ref{relations}) of the algebra
$O_q(M(n))$, it is easy to show the following lemma.

\begin{Lem} (1) The mapping
\begin{eqnarray}^-:Z_{ij}&\mapsto &Z_{ij}\\\nonumber
q&\mapsto &q^{-1}\end{eqnarray}
extends to an algebra anti-automorphism of the algebra
$O_{q}(M(n))$ as an algebra over ${\mathbb Q}$.

(2) The mapping
\begin{equation}\sigma: Z_{ij}\mapsto Z_{ji}\end{equation}
extends to an algebra automorphism of the algebra $O_{q}(M(n))$ as an algebra over $K={\mathbb Q}(q)$.
\end{Lem}

For any $A=(a_{ij})_{n\times n}\in M_n({\mathbb Z}_+)$.
Let
$$ro(A)=(\sum_j a_{1j},\cdots,\sum_ja_{nj})=(r_1,r_2,\cdots,r_n)$$
which is called the row sum of $A$
 and
$$co(A)=(\sum_j a_{j1},\cdots,\sum_j a_{jn})=(c_1,c_2,\cdots,c_n)$$
which is called the column sum of $A$.

For any matrix $A=(a_{ij})_{i,j=1}^n\in M_n({\mathbb Z}_+)$, a
monomial having the factors of $Z^A$ in arbitrary order. Then its
expansion in terms of monomials $Z^B$ only involves terms where
$ro(B)=ro(A)$ and $co(B)=co(A)$. Let $Pr(A, s,t)=\sum_{i\le s,
j\le t}a_{ij}$. Then $Pr(A,s,t)\ge Pr(B,s,t)$ for any $s,t\le n$
and matrix $B$ appeared in the expansion considered above.

From the defining relations (\ref{relations}) of the algebra
$O_q(M(n))$, we have

\begin{equation}\label{bar}\overline{Z^A}=E(A)Z^A+\sum_{B}c_B(A)Z^B,
\end{equation}
where
$$E(A)=q^{-2(\sum_i\sum_{j>k}a_{ij}a_{ik}+\sum_i\sum_{j>k}a_{ji}a_{ki})}$$
and $B<A$, $ro(B)=ro(A)$, $co(B)=co(A)$, $c_B(A)\in{\mathbb Z}[q,q^{-1}]$,
$\le$ is the lexicographic ordering.

For a pair of vectors $R,C\in{\mathbb Z}_+^n$, denote by $M(R,C)$ the subspace of
$ O_q(M(n))$ spanned by $Z^A$ with $ro(A)=R,$ and $co(A)=C$. Note that $M(R,C)$ is
$^-$ invariant and
$ O$$_q(M(n))=\oplus_{R,C} M(R,C)$.

Let $D(A)=q^{-\sum_i \sum_{j>k}a_{ij}a_{ik}-\sum_i
\sum_{j>k}a_{ji}a_{ki}}$ and let $Z(A)=D(A)Z^A.$ Set
 $$L^*=\oplus_{A\in M_n({\mathbb Z}_+)}{\mathbb Z}[q]Z(A).$$

\begin{Thm}\label{basis}
There is
a unique basis
$B^*=\{b(A)|A\in M_n({\mathbb Z}_+)\}$ of $ L^*$ determined by the following conditions:
\begin{enumerate}
\item $\overline{b(A)}=b(A)$ for all $A$.

\item $b(A)=Z(A)+\sum_{B<A} h_B(A)Z(B)$ where $h_B(A)\in q{\mathbb Z}$$[q]$ and $ro(B)=ro(A), co(B)=co(A)$.
\end{enumerate}
\end{Thm}

\pf    We rewrite the equation (\ref{bar}) in terms of $Z(A)$, then
\begin{equation}\overline{Z(A)}=\sum_Ba_{AB}Z(B),\end{equation}
where $a_{AA}=1$, $a_{AB}\in{\mathbb Z}[q,q^{-1}]$ and $a_{AB}=0$ unless
$B\le A$, where $\le$ is the lexicographic ordering. By Theorem 1.2 of \cite{d1}, there
is an IC-basis with respect to the triple $(\{Z^A\mid A\in M_n({\mathbb Z}_+)\},
^-, \le)$ determined by the relation stated in the context of the theorem.
\qed

The quantum determinant ${\det}_q$ is defined as follows:

\begin{equation} {\det}_q=
\Sigma_{\sigma\in S_n}(-q^2)^{l(\sigma)}Z_{1\sigma(1)}Z_{2\sigma(2)} \cdots Z_{n\sigma(n)}.\end{equation}
It is known that $det_q$ is a central element of the algebra  $ O_q(M(n))$.

For later reference we now introduce some terminology. Let $m\le n$ be a positive integer.
 Given any two subsets $I=\{i_1,i_2,\cdots,i_m\}$ and $J=\{j_1,j_2,\cdots,j_m\}$ of
 $\{1,2,\cdots,n\}$, each having cardinality $m$, it is  clear that the subalgebra of
 $ O$$_q(M(n))$ generated by the elements $Z_{i_rj_s}$ with $r,s=1,2,\cdots,m$,
 is isomorphic to $ O$$_q(M(m))$, so we can talk about its determinant.
 Such a determinant is called a quantum minor, and will be denoted by ${\det}_q(I,J)$.

Let $I, J$ be two subsets of  $\{1,2,\cdots,n\}$ with the same cardinality.
Obviously, the dual canonical basis of the subalgebra generated by $Z_{ij}$
for $i\in I$ $j\in J$ is a subset of the basis  $B^*$ of the algebra $O_q(M(n))$.
More generally, If $(u,v)\le (s,t)$, then the subalgebra $O_q(M(n))^{(u,v)}_{(s,t)}$
generated by $Z_{i,j}$, for $(u,v)\le (i,j)\le (s,t)$, is $^-$ invariant and one can
construct a basis analogous to the construction of the  basis considered in
Theorem \ref{basis}, and obviously the resulting basis of $O_q(M(n))^{(u,v)}_{(s,t)}$
is a subset of the basis $B^*$.

\begin{Lem} The quantum determinant $det_q$ is an element of the basis $B^*$. Furthermore, any quantum minor is also an element of the dual canonical basis.\end{Lem}

\pf   We only need to show that $det_q$ is $^-$ invariant. It is well known that the center of the algebra $O_q(M(n))$ is generated by the quantum determinant \cite{nmy}. Note that
\begin{eqnarray}\overline{det_q}Z_{ij}&=&\overline{\overline{Z_{ij}}det_q}=\overline{det_qZ_{ij}}\\\nonumber
&=&Z_{ij}\overline{det_q},\end{eqnarray}
for any $i,j$. Hence, $\overline{det_q}$ is a polynomial of $det_q$. Therefore,
$$\overline{det_q}=det_q$$
by comparing the leading terms.\
\qed

\begin{Cor}\label{transpose}The  basis $B^*$ is $\sigma$ invariant. More precisely,

$$\sigma(b(A))=b(A^T),$$
for all $A\in M_n({\mathbb Z}_+)$, where $A^T$ is the transposition of $A$.\end{Cor}

\pf
Let $ b(A)$ be an element of the dual canonical basis $B^*$ of the form given in Theorem \ref{basis} (2). Then it follows that
all of the matrices $B$ appearing in the expansion of $b(A)$ are obtained
from $A$ by a sequence of  $2\times 2$ submatrix transformations of the following form:

\begin{equation}\begin{pmatrix}a_{ij}&a_{it}\\a_{sj}&a_{st}\end{pmatrix}
\longrightarrow
\begin{pmatrix}a_{ij}-1&a_{it}+1\\a_{sj}+1&a_{st}-1\end{pmatrix},\end{equation}
if both $a_{ij}$ and $a_{st}$ are positive. Hence $B^T$ can be
obtained from $A^T$ by  a sequence of the  submatrix
transformations of the form:

\begin{equation}\begin{pmatrix}a_{ji}&a_{ti}\\a_{js}&a_{ts}\end{pmatrix}
\longrightarrow
\begin{pmatrix}a_{ji}-1&a_{ti}+1\\a_{js}+1&a_{ts}-1\end{pmatrix}.\end{equation}
Especially, $B^T\le A^T$.  Note that the monomials $Z^{B^T}$ and $\sigma(Z^B)$ have the same factors but could be in different order. However, two generators $Z_{ij}$ and $Z_{st}$ appear in the monomials but in different order must satisfy the third relation in (\ref{relations}).  Hence, $Z^{B^T}=\sigma(Z^B)$
\begin{equation}\sigma(b(A))=Z(A^T)+\sum_{B}h_B(A)Z(B^T)\end{equation}
with $h_B(A)\in q{\mathbb Z}[q]$. Clearly,

$$\overline{\sigma(b(A)}=\sigma(b(A))$$
since $\sigma$ and $^-$ commute with each other. \qed

Denote by $I_n$ the $n\times n$ identity matrix.

\begin{Lem}For any $A\in M_n({\mathbb Z}_+)$,
$$Z(A)det_q=Z(A+I_n)\quad \bmod q L^*.$$
\end{Lem}

\pf   For $i<s, j<t$, we have

$$Z_{st}^mZ_{ij}=Z_{ij}Z_{st}^m+(q^{2-4m}-q^2)Z_{it}Z_{sj}Z_{st}^{m-1}.$$
Recall that
$${\det}_q=
\Sigma_{\sigma\in S_n}(-1)^{l(\sigma)}q^{2l(\sigma)}
Z_{1\sigma(1)}Z_{2\sigma(2)} \cdots Z_{n\sigma(n)}.$$
When we compute $Z(A)det_q$, we only have  to deal with those coefficients of
the form $q^{-2a}$ with $a$ a positive integer. Assume that
$$Z(A)det_q=\sum a_BZ(B).$$
Clearly, $a_B\in{\mathbb Z}[q,q^{-1}]$ and the leading term is
$Z(A+I_n)$ and those matrix $B$ appeared in the expression has at
least one nonzero entry in each row and each column. We need to
compute $Z(A)(-1)^{l(\sigma)}q^{2l(\sigma)}
Z_{1\sigma(1)}Z_{2\sigma(2)} \cdots Z_{n\sigma(n)},$ for all
$\sigma\in S_n$. From the expression of the quantum determinant we
see that there are four possibilities to produce coefficients of
the form $q^{-2a}$ with $a$ a positive integer.

Case 1. $Z_{st}^mZ_{sj}=q^{-2m}Z_{sj}Z_{st}^m$ where $t>j$ but no $Z_{it}$ behind. Then $q^{2m}$ will be absorbed by $D(B)$ where $Z(B)$ is the resulted term.

Case 2. $Z_{st}^mZ_{it}=q^{-2m}Z_{it}Z_{st}^m$ where $s>i$ but no $Z_{sj}$ appeared before. Then $q^{2m}$ will be absorbed by $D(B)$ where $Z(B)$ is the resulted term.

Case 3. Both $Z_{st}^mZ_{sj}=q^{-2m}Z_{sj}Z_{st}^m$ where $t>j$ and
$Z_{st}^mZ_{it}=q^{-2m}Z_{it}Z_{st}^m$ where $s>i$ happened. Then we get
$q^{-4m}$. However, we will see that it will be cancelled by a term in the next case. To this end, we need to remember that the terms we are dealing with are from $Z(A) Z_{1\sigma(1)}Z_{2\sigma(2)} \cdots Z_{n\sigma(n)}.$ Note that $l(\sigma (jt))=l(\sigma)-1$.

Case 4. $Z_{st}^mZ_{ij}= Z_{ij}Z_{st}^m+(q^{2-4m}-q^2)Z_{it}Z_{sj}Z_{st}^{m-1}$
where $s>i, t>j$. Then the coefficient $q^{2-4m}$ will be cancelled by a term in case 3.

Hence, the coefficients $a_B$ are all in $q{\mathbb Z}[q]$ except $a_A$ which is $1$.
\qed

The following proposition follows directly from the above lemma.

\begin{Prop}The basis $B^*$ is invariant under the multiplication of $det_q$. More precisely,
$$b(A)det_q=b(A+I_n)$$
for all $A\in M_n({\mathbb Z}_+)$.\end{Prop}

By using this proposition, we can determine $b(A)$, if $A$ is a diagonal matrix. Let $A=diag(a_1, a_2,\cdots,a_n)$. We may assume that $a_1\le a_2\le \cdots \le a_n$ without loss  of generality.  Then

$$b(A)=\Pi_{i=1}^n det_{q, i}^{a_i-a_{i-1}}$$
where $det_{q,i}$ is the quantum determinant of the subalgebra
generated by $Z_{st}$ for $s,t=i,\cdots,n$, and where we put $a_0=0$.

\medskip

\section{some subalgebras}
\medskip

In this section, we study the multiplicative property of the basis
$B^*$. Similar to the proof of Proposition 2.9, we get

\begin{Lem}Let $b_1, b_2\in B^*$. If $b_1b_2\sim b$ for some $b\in B^*$,
 then $b_1b_2\sim b_2b_1$.\end{Lem}

Divide the matrix  by a broken-line $\xi$ which  consists lines
determined by the equations $ax+by=m$ for $a,b\in {\mathbb Z}_+$
and $m\in\mathbb N$ (each line has non-positive slope). Recall
that $I(n)=\{(i,j)\mid i,j=1,\dots,n\}$. Let
$$I_1=\{(x,y)\in I(n)\mid (x,y)\text{ is in the left upper side of the broken line } \xi\},$$
and let $I_2$ be the complement of $I_1$ in $I(n)$.

Let $O_i$ be the subalgebra of $O_q(M(n))$ generated by $Z_{xy}$
for $(x,y)\in I_i$, $i=1,2$. One can easily see that $O_i$ is
determined by the generators $Z_{xy}$ and relations
(\ref{relations}). Hence, the algebra $O_i$ is closed under the
bar action and therefore there is a basis $B_i^*$ of the
sub-lattice $L_i^*$ of the ${\mathbb Z}[q]$-lattice $L^*$ spanned
by $\{Z(A)\mid A=(a_{xy})\in M_n({\mathbb Z}_+), a_{xy}=0 \text{
if }(x,y)\in I_{3-i}\}$. Clearly, $B^*_i$ is a subset of $B^*$
consists of those $b(A)$ for $A=(a_{xy})\in M_n({\mathbb Z}_+),
a_{xy}=0 \text{ if }(x,y)\in I_{3-i}$.

Write
$$A=A^+ +A^-,$$
Where the entries of $A^+$ in the left upper side of the broken
line $\xi$ are zero and the entries of $A^-$ in the right lower
side (including the broken line $\xi$) are zero. Then

\begin{Thm}$b(A)\sim b(A^+)b(A^-)$ if and only if $ b(A^+)b(A^-)\sim b(A^-)b(A^+)$
\end{Thm}

\pf    If $b(A)=q^ab(A^+)b(A^-)$ for some integer $a$, then
$ b(A^-)b(A^+)\sim  b(A^+)b(A^-)$ by the above lemma.

For
$$b(A^+)=Z(A^+)+\sum_{B^+}a_{B^+A^+}Z(B^+), $$
and
$$b(A^-)=Z(A^-)+\sum_{B^-}a_{B^-A^-}Z(B^-), $$
where $ a_{B^+A^+}, a_{B^-A^-}\in q{\mathbb Z}[q]$. Assume that
$b(A^+)b(A^-)=q^a b(A^-)b(A^+)$, for some integer $a$ which can be
computed by only considering the leading terms. From the defining
relations (\ref{relations}), the integer $a$ must be even, say,
$a=2m$. Then $q^{-m} b(A^+)b(A^-)$ is bar-invariant with leading
term $Z(A)$. Note that the coefficients we encounter only depend
on the row sums and column sums. Actually,
$m=\sum_{j}(r_j^+r_j^-+c_j^+c_j^-)$ where $(r_1^+,\cdots, r_n^+)$
and $(c_1^+,\cdots,c_n^+)$ (resp. $(r_1^-,\cdots, r_n^-)$ and
$(c_1^-,\cdots,c_n^-)$ are the row sum and column sum respectively
of $A^+$ (resp of $A^-$). Then all term produce the same $m$.
Therefore,
$$b(A)=q^{-m}b(A^+)b(A^-)$$
by  Theorem \ref{basis}.

\qed

\medskip

\section{some quantum minors}

\medskip

Let ${\det}_q(t)={\det}_q(\{1,\cdots,t\},\{n-t+1,\cdots,n\})$, for $t=1,2,\cdots,n$.
\smallskip

 Let $M_t^-=\{(i,j)\in \mathbb N$$^2\mid 1\le i\le t\text{ and }1\le j\le n-t\}$,
 $M_t^+=\{(i,j)\in \mathbb N$$^2\mid t+1\le i\le n\text{ and }n-t+1\le j\le n\}$,
 $M_t^l=\{(i,j)\in \mathbb N$$^2\mid t+1\le i\le n\text{ and }1\le j\le n-t \}$,
 and $M_t^r=\{(i,j)\in \mathbb N$$^2\mid 1\le i\le t\text{ and }n-t+1\le j\le n\}$.
The following result was proved in \cite{jz}.

\begin{Lem} For any $i, j, t$,
\begin{eqnarray}
Z_{ij}{\det}_q(t)&=&{\det}_q(t)Z_{ij} \text{ if }(i,j)\in M_t^l\cup M^r_t,
\\\nonumber Z_{ij}{\det}_q(t)&=&q^2{\det}_q(t)Z_{ij}
\text{ if }  (i,j)\in M_t^-,\text{ and }\\\nonumber
Z_{ij}{\det}_q(t)&=&q^{-2}{\det}_q(t)Z_{ij} \text{ if }(i,j)\in M^+_t. \end{eqnarray}
\end{Lem}

Let $E_t=\begin{pmatrix}0&I_t\\0&0\end{pmatrix}$ and let
$ q^{{\mathbb Z}}$$B^*=\{q^ab(A)|\text{ for all } A \text{ and }a\in{\mathbb Z}\}$.

\begin{Thm} The set $q^{{\mathbb Z}}$$B^*$ is invariant under the multiplication of the
quantum minors $det_q(t)$ and $\sigma(det_q(t))$. More precisely,

\begin{equation}\label{c1}b(A)det_q(t)=q^{r_1+\cdots +r_t-c_{n-t+1}-\cdots c_n}b(A+E_t).\end{equation}
\begin{equation}\label{c2}b(A)\sigma(det_q(t))=q^{c_1+c_2+\cdots +c_t-r_{n-t+1}-\cdots -r_n}b(A+E^T_t),\end{equation}
where $E^T_t$ is the transposition of the matrix $E_t$.
\end{Thm}

\medskip

\pf

For any $A\in M_n({\mathbb Z}_+)$, write
$$A=\begin{pmatrix}A_{11}&A_{12}\\A_{21}&A_{22}\end{pmatrix},$$
where $A_{11}$ is a $t\times (n-t)$ submatrix, $A_{12}$ is a $t\times t$ submatrix, $A_{21}$ is a $(n-t)\times (n-t)$ submatrix and $A_{22}$ is a $(n-t)\times t$ submatrix. Then
two monomials $Z^A$ and $Z^{A_{11}}Z^{A_{12}}Z^{A_{21}}Z^{A_{22}}$ have the same factors but could be in different order. However, the first monomial can be obtained from the second one by applying the third defining relation of the algebra $O_q(M(n))$. Hence,
$$Z^A=Z^{A_{11}}Z^{A_{12}}Z^{A_{21}}Z^{A_{22}},$$
and
$$Z(A)det_q(t)=
q^{-2\sum_{i\ge t+1, j\ge n-t+1}a_{ij}}D(A)Z^{A_{11}}
Z^{A_{12}}det_q(t)Z^{A_{21}}Z^{A_{22}}.$$
Apply the above lemma to $Z^{A_{12}}det_q(t)$, then we get
$$ Z(A_{12})det_q(t)=Z(A_{12}+I_t)+\sum_{B_{12}}h_{B_{12}}(A_{12})Z(B_{12}),$$
where $B_{12}\le A_{12}+I_t$,  $B_{12}$ and $A_{12}+I_t$ have the same row sums and column sums. Hence
$$ Z(A)det_q(t)=
q^{-2\sum_{i\ge t+1, j\ge n-t+1}a_{ij}}D(A)
D(A_{12})^{-1}$$
$$D(A_{12}+I_t)D(A+E_t)^{-1}
Z(A+E_t)$$
$$+\sum_{B_{12}}h_{B_{12}}(A_{12})
q^{-2\sum_{i\ge t+1, j\ge n-t+1}a_{ij}}D(A)D(A_{12})^{-1}$$
$$D(B_{12})D(\begin{pmatrix}A_{11}&B_{12}\\A_{21}&A_{22}\end{pmatrix})^{-1}
Z(\begin{pmatrix}A_{11}&B_{12}\\A_{21}&A_{22}\end{pmatrix}).$$

By direct computation, one can show that the dependence of

$$D(B_{12})D(\begin{pmatrix}A_{11}&B_{12}\\A_{21}&A_{22}\end{pmatrix})^{-1}$$
on the matrix entries of $B_{12}$ is only a dependence on the row and column sums.

Then one deduces that
$$Z(A)det_q(t)=q^{r_1+\cdots+r_t-c_{n-t+1}-\cdots-c_n}(Z(A+E_t)+
\sum_{D,D<A+E_t}c(D, A)Z(D))$$
with $c(D, A)\in q{\mathbb Z}$$[q]$.

For a basis element $b(A)$ of the form given in
Theorem\ref{basis}, we then deduce that

\begin{eqnarray} b(A)det_q(t)& = & q^{r_1+\cdots+r_t-c_{n-t+1}-\cdots-c_n}\\\nonumber
& &[(Z(A+E_t)+\sum_{D, D<A+E_t} c_D(A)Z(D)) \\\nonumber
& &+\sum_{B, B<A}h_B(A)(Z(B+E_t)+\sum_{D,D<B+E_t}c_D(B)Z(D))]\end{eqnarray}
with $c_D(A), c_D(B)\in q{\mathbb Z}$$[q]$. By
$$b(A)det_q(t)= q^{2(r_1+\cdots+r_t-c_{n-t+1}-\cdots-c_n)}det_q(t)b(A),$$
we see that $(Z(A+E_t)+\sum_{D, D<A+E_t} c_D(A)Z(D))+\sum_{B,
B<A}h_B(A)(Z(B+E_t)+ \sum_{D,D<B+E_t}c_D(B)Z(D))$ is $^-$
invariant, and it must be the basis element $b(A+E_t)$. Finally,
apply the algebra automorphism $\sigma$. Then the second statement
follows from Corollary \ref{transpose}\qed

\begin{Cor} Let $A=\begin{pmatrix}a_{1}&b_2&b_3&\cdots &b_{n-1}&b_{n}\\c_2&a_{2}&b_2&\cdots &b_{n-2}&b_{n-1}\\
c_3&c_2&a_{3}&\cdots &b_{n-3}&b_{n-2}\\ \quad\quad\cdots&\cdots&\cdots&\cdots\\
c_n&c_{n-1}&c_{n-2}&\cdots&c_2&a_{n}\end{pmatrix}$. Then the basis element $$b(A)\sim \Pi_{t=1}^{n-1}det_q(t)^{b_{n-t+1}}\sigma(det_q(t))^{c_{n-t+1}} b(diag(a_{1},a_{2},\cdots,a_{n})).$$ \end{Cor}

\pf    Successively peel off the off-diagonals of $A$ by
(\ref{c1}) and (\ref{c2}). \qed

 \begin{Def} A matrix $A=(a_{ij})\in M_n({\mathbb Z}_+)$ is called a
 ladder
 if $a_{ij}\ge a_{i+1, j+1}$ for all $i,j$. \end{Def}

 Let $A$ be a ladder. Successively peel off the off-diagonals of $A$ by (\ref{c1}) and
 (\ref{c2}), the basis element $b(A)$ is equivalent to a product of the
 quantum minors $det_q(t)$ and $\sigma(det_q(t))$ and a basis
 element $b(\begin{pmatrix}A_{n-1}&0\\0&0\end{pmatrix})$, where
 $A_{n-1}$ is a ladder of size $n-1$. Repeatedly, the basis
 element $b(A)$ can be written as a product of some quantum minors
 which are $q$-commuting with each other.

\medskip

\section{ the coincidence of two bases}

\medskip

Let $g$ be the finite dimensional simple Lie algebra of type
$A_{n-1}$ and let $\Lambda_1,\cdots,\Lambda_{n-1}$ be the
fundamental dominant weights. For any dominant weight $\lambda$,
the irreducible highest weight module $L(\lambda)$ occurs as a
sub-quotient of a suitable power of the natural representation
$L(\Lambda_1)$. The simple modules $L(\Lambda_1)$ and
$L(\Lambda_{n-1})$ are dual to each other and are of dimension
$n$. Let $e_1,e_2,\cdots,e_n$ be the standard basis of
$L(\Lambda_1)$ and let $e_1^*,e_2^*,\cdots,e_n^*$ be the dual
basis of $L(\Lambda_{n-1})$. Then it is well-known that the matrix
coefficients $X_{ij}=e_i^*\otimes e_j$ satisfy the following
relations:

\begin{eqnarray}X_{ij}X_{ik}&=&q^2X_{ik}X_{ij} \text{ if } j<k,\\
 X_{ij}X_{kj}&=&q^2X_{kj}X_{ij} \text{ if }i<k,\\
X_{ij}X_{st}&=&X_{st}X_{ij}\text{ if } i>s, j<t,\\
X_{ij}X_{st}&=&X_{st}X_{ij}+(q^2-q^{-2})X_{it}X_{sj} \text{ if }
i<s, j<t,\\\nonumber \Sigma_{\sigma\in
S_n}&(-q^2)^{l(\sigma)}&X_{1\sigma(1)}X_{2\sigma(2)} \cdots
X_{n\sigma(n)}=1.
\end{eqnarray}

Since the basis $B^*$ is invariant under the multiplication of the
quantum determinant,  we gets a basis $K^*$ of
$O_q(SL_n)$($=A_q(g)$), by setting the quantum determinant to one
. Clearly, the anti-automorphism $-$ induces the anti-automorphism
$\phi$ of $O_q(SL(n))$ (see Lemma 2.9). Let $X(A)$ be the image of
$Z(A)$ in $O_q(SL(n))$. Then
$$\{X(A)|trA=0\}$$
is a basis of $O_q(SL(n))$.

\begin{Lem} The matrix coefficients $X_{ij}$ are both invariant
under $-$ (the bar action of $A_q(g)$) and $\phi$.\end{Lem}

\pf It is known that $\{e_1,e_2,\cdots,e_n\}$ (resp.
$\{e_1^*,e_2^*,\cdots,e_n^*\}$) is the canonical basis of
$L(\Lambda_1)$ (resp.  of $L(\Lambda_{n-1})$). Therefore, $e_i$
and $e_j^*$ are invariant under the bar action of $L(\Lambda_1)$
and $L(\Lambda_{n-1})$ respectively. Hence, the matrix
coefficients $X_{ij}$ are $-$ invariant. Note that $\Lambda_1$ and
$\Lambda_{n-1}$ are minuscule dominant weights so the left weight
$\lambda_l$ (resp. the right weight $\lambda_r$) of $X_{ij}$ is
conjugate to $\Lambda_1$ (resp. $\Lambda_{n-1}$) under the action
of the Weyl group which implies that
$(\lambda_l,\lambda_l)-(\lambda_r,\lambda_r)=(\Lambda_1,\Lambda_1)
-(\Lambda_{n-1},\Lambda_{n-1})=0$. \qed

The basis $K^*$ can be described similarly to the theorem
\ref{basis} by replacing $Z_{ij}$ by $X_{ij}$ and $-$ by $\phi$.

\begin{Thm}
There is a unique basis $\tilde{B^*}=\{\tilde{b(A)}|A\in
M_n({\mathbb Z}_+), trA=0\}$ of $ \tilde{L^*}=\oplus_A{\mathbb
Z}[q]X(A)$ determined by the following conditions:
\begin{enumerate}
\item $\phi\tilde{b(A)}=\tilde{b(A)}$ for all $A$.

\item $\tilde{b(A)}=X(A)+\sum_{B<A} h_B(A)X(B)$ where $h_B(A)\in
q{\mathbb Z}$$[q]$ and $ro(B)=ro(A), co(B)=co(A)$.
\end{enumerate}
\end{Thm}

Let ${\mathbb R}^n$ be the $n$ dimensional Euclidean space with
standard orthogonal  basis $\epsilon_1, \epsilon_2, \cdots,
\epsilon_n$. It is well-known that the root system of type
$A_{n-1}$ is a subset of ${\mathbb R}^n$ with simple roots
$\alpha_i=\epsilon_i-\epsilon_{i+1}$, for $i=1,2,\cdots, n-1$.

The $U_q(g)$ bi-module structure can ba written explicitly (see
also \cite{nym}).

For homogeneous elements $x,$ and $y$ with weights $(\lambda_l,
\lambda_r)$ and $(\mu_l, \mu_r)$ respectively.

The left action is defined by

$$E_iX_{st}=\delta_{is}X_{s-1, t}, F_iX_{st}=\delta_{i,
s+1}X_{s+1, t}, K_iX_{st}=q^{2(\epsilon_s, \alpha_i)}X_{st}$$

with Leibniz rule

$$E_i(xy)=E_i(x)y+q^{2(\lambda_l, \alpha_i)}xE_i(y),$$

$$F_i(xy)=xF_i(y)+q^{-2(\mu_l, \alpha_i)}F_i(x)y,$$

$$K_i(xy)=q^{2(\lambda_l+\mu_l, \alpha_i)}xy,$$

The right action is defined by

$$X_{st}E_i=\delta_{i,s+1}X_{s+1, t}, X_{st}F_i=\delta_{i,
s}X_{s-1, t}, X_{st}K_i=q^{2(\epsilon_s, \alpha_i)}X_{st}$$

with Leibniz rule

$$(xy)E_i=(x)E_i y+q^{2(\lambda_r, \alpha_i)}x(y)E_i,$$

$$(xy)F_i=x(y)F_i+q^{-2(\mu_r, \alpha_i)}(x)F_i y,$$

$$(xy)K_i=q^{2(\lambda_r+\mu_r, \alpha_i)}xy,$$

Denote by the same notation the image of $det_q(i)$ in
$O_q(SL(n))$.

For $\lambda=m_1\Lambda_1+m_2\Lambda_2+\cdots
+m_{n-1}\Lambda_{n-1}$, where $\Lambda_1,\Lambda_2,\cdots,
\Lambda_{n-1}$ are fundamental weights.  The module
$L(\lambda)\otimes L^*(\lambda)$ is cyclic on $v_\lambda\otimes
v_\lambda^*$. The lattice ${\cal L}(\lambda)\otimes {\cal
L}^*(\lambda)$ is generated by  $v_\lambda\otimes v_\lambda^*$
(under the actio of upper Kashiwara operators) which corresponds
to
$$\Pi_idet_q(i)^{m_i}$$ which is an element in the basis $K^*$.
The quantum minor $det_q(t)$ is annihilated by the left action of
$E_i$ for $i=1,2,\cdots, t-1$ and the right action of $F_j$ for
$j=n-1, n-2, \cdots, n-t+1$.

By the action of the upper kashiwara operators, the normalized
monomial $X(A)$ and the element
$q^{\frac{1}{2}((\lambda_l,\lambda_l)-(\lambda_r,\lambda_r))}b^\prime$
are in the same ${\mathbb Z}[q]$-lattice. By the uniqueness of
Lusztig's construction. The bases $K^*$ and $L^*$ are coincide.

\bibliographystyle{amsplain}

\end{document}